\documentclass[11pt]{article}
\usepackage{amsfonts}
\usepackage[dvips]{graphics}
\usepackage{epsfig}
\usepackage{graphicx}
\usepackage{psfrag}
\usepackage{amssymb,subfigure}
\usepackage{amsthm}
\usepackage{mathrsfs}
\usepackage{color}
\usepackage{algorithm2e}

\def\R{\mathbb R}

\def\tmu{\tilde{\mu}}
\def\tg{\tilde{\gamma}}

\def\ua{ u}
\def\ub{\tilde{u}}
\def\uc{\tilde{u}}

\def\epsilon{\varepsilon}
\def\ds{\displaystyle}

\newcommand{\be}{\begin{equation}}
\newcommand{\ee}{\end{equation}}
\newcommand{\baa}{\begin{array}}
\newcommand{\eaa}{\end{array}}
\newcommand{\ba}{\begin{eqnarray}}
\newcommand{\ea}{\end{eqnarray}}

\newtheorem{lemma}{Lemma}[section]
\newtheorem{theorem}[lemma]{Theorem}
\newtheorem{corollary}[lemma]{Corollary}

\newtheorem{proposition}[lemma]{Proposition}
\newtheorem{remark}[lemma]{Remark}
\newtheorem{remarks}[lemma]{Remarks}

\title{On the determination of the nonlinearity from localized measurements in a reaction-diffusion equation}

\author{L. Roques$^{\hbox{\small{ a,*}}}$ and M. Cristofol$^{\hbox{\small{ b }}}$
\\
\footnotesize{$^{\hbox{a }}$UR 546 Biostatistique et Processus Spatiaux, INRA, F-84000
Avignon, France}\\
\footnotesize{$^{\hbox{b }}$Aix-Marseille Universit\'e, LATP, F-13397 Marseille, France} \\
\footnotesize{$^{\hbox{*} }$ Author for correspondence. (lionel.roques@avignon.inra.fr)}
}

\date{}

\begin{document}

\maketitle

\begin{abstract}
This paper is devoted to the analysis of some uniqueness properties of a classical reaction-diffusion equation of  Fisher-KPP type, coming from  population dynamics in heterogeneous environments. We work in a one-dimensional interval $(a,b)$ and we assume a nonlinear term of the form $u \, (\mu(x)-\gamma u)$ where $\mu$ belongs to a fixed subset of $C^{0}([a,b])$. We prove that the knowledge of $u$ at $t=0$ and of $u$, $u_x$  at a single point $x_0$ and for small times $t\in (0,\varepsilon)$ is sufficient to completely determine the couple $(u(t,x),\mu(x))$ provided $\gamma$ is known. Additionally, if $u_{xx}(t,x_0)$ is also measured for $t\in (0,\varepsilon)$, the triplet $(u(t,x),\mu(x),\gamma)$ is also completely determined. Those analytical results are completed with numerical simulations which show that, in practice, measurements of $u$ and $u_x$ at a single point $x_0$ (and for $t\in (0,\varepsilon)$) are sufficient to obtain a good approximation of the coefficient $\mu(x).$ These numerical simulations also show that the measurement of the derivative $u_x$ is essential in order to accurately determine $\mu(x)$.
\end{abstract}


%
%%%%%%%%%%%%%%%%%%%%%%%%%%%%%%%%%%%%%%%%%%%%%%%%%%%%%%%%%%%%
\section{Introduction and ecological background}
%%%%%%%%%%%%%%%%%%%%%%%%%%%%%%%%%%%%%%%%%%%%%%%%%%%%%%%%%%%%
%

Reaction-diffusion models (hereafter RD models), although they
sometimes bear on simplistic assumptions such as infinite velocity
assumption and completely random motion of animals
\cite{Hol93}, are not in disagreement with certain dispersal
properties of populations observed in natural as well as
experimental ecological systems, at least qualitatively
\cite{ShiKaw97,Tur98,Mur02,OkuLev02}. In fact, since the work of Skellam \cite{Ske51}, RD theory has been the main analytical framework to study spatial spread of biological organisms, partly because it benefits from a well-developed mathematical theory.

The idea of modeling population dynamics with such models has emerged at the beginning of the 20$^{\hbox{th}}$ century, with random walk theories of organisms, introduced by  Pearson and Blakeman \cite{PeaBla06}. Then, Fisher  \cite{Fis37} and   Kolmogorov, Petrovsky, Piskunov \cite{kpp} independently used a reaction-diffusion equation as a model for population genetics.
The corresponding equation is
\begin{equation}\label{eqkpp}
\frac{\partial u}{\partial t}-D  \frac{\partial^2 u}{\partial x^2}=u\, (\mu- \gamma u), \ t>0, \ x\in(a,b)\subset\R,
\end{equation}
where $u=u(t,x)$ is the population density at time $t$ and space
position $x$, $D$ is the diffusion coefficient, and  $\mu$ and
$\gamma$ respectively correspond to the \emph{constant} intrinsic growth
rate and intraspecific competition coefficients. In the 80's, this
model has been extended to heterogeneous environments by Shigesada
et al. \cite{ShiKaw86}. The corresponding model  is of the type: \be\label{eqskt_gale}
\frac{\partial u}{\partial t}- \frac{\partial}{\partial x}\left(D(x)\frac{\partial u}{\partial x}\right)=u\, (\mu(x)-\gamma(x) u), \ t>0, \ x\in(a,b).\ee
The coefficients $\mu(x)$
and $\gamma(x)$ now depend on the space variable $x$ and can
therefore include some effects of environmental heterogeneity.
More recently, this model revealed that the heterogeneous
character of the environment played an essential role on species
persistence and spreading, in the sense that for different spatial
configurations of the environment, a population can survive or
become extinct and spread at different speeds, depending on the habitat spatial structure
(\cite{ShiKaw97}, \cite{BerHamRoq05a}, \cite{CanCos03},\cite{ElsHamRoq09}, \cite{HamFay10} ,\cite{RoqHam07}, \cite{RoqSto07}).
Thus, determining the coefficients in model (\ref{eqskt_gale})  is an important question, even for areas other than ecology (see \cite{Xin00} and references therein).

In this paper, we focus on the case of constant coefficients $D$ and $\gamma$:
\be\label{eqskt}
\frac{\partial u}{\partial t}-D  \frac{\partial^2 u}{\partial x^2}=u\, (\mu(x)-\gamma u), \ t>0, \ x\in(a,b),
\ee
 and
we address the question of the uniqueness of couples $(u,\mu(x))$ and triples $(u,\mu(x),\gamma)$ satisfying (\ref{eqskt_gale}),
given a localized measurement of $u$.

Uniqueness results of this type have been obtained for reaction-diffusion models, through the
Lipschtiz stability of the coefficient with respect to the solution $u$. Lipschtiz  stability is generally obtained
by using the method of Carleman estimates \cite{BukKli81}. Several publications starting from the paper by Isakov \cite{Isa93} and including the recent overview of the method of Carleman estimates applied to inverse coefficients problems \cite{KliTim04} provide results for the case of multiple measurements. The particular problem of the uniqueness of the couple $(u,\mu(x))$  satisfying (\ref{eqskt}) given  such multiple measurements has been investigated, together with Lipschtiz stability, in a previous work \cite{CriRoq08}. Placing ourselves in a bounded domain $\Omega$ of $\R^N$ with Dirichlet boundary conditions, we had to use the following measurements: (i) the density $u(0,x)$ in $\Omega$ at $t=0$;  (ii)  the density $u(t,x)$ for $(t,x)\in (t_0,t_1)\times \omega$, for some times $0<t_0<t_1$ and a subset $\omega \subset \subset \Omega$; (iii) the density $u(\theta,x)$  for all $x\in \Omega$, at some time $\theta\in (t_0,t_1)$.

Although the  result of \cite{CriRoq08} allows to determine $\mu(x)$ using partial measurements of $u(t,x),$ assumption (iii) implies that $u$ has to be known in the whole set $\Omega$. This last measurement (iii) is a key assumption in several other papers  on uniqueness and stability of solutions to     parabolic equations with respect to parameters (see Imanuvilov and Yamamoto \cite{ImmYam98}, Yamamoto and Zou \cite{YamZou01}, Belassoued and Yamamoto \cite{BelYam06} for scalar equations  and   Cristofol, Gaitan and Ramoul \cite{CriGai06} or Benabdallah, Cristofol, Gaitan and Yamamoto \cite{BenCri08} for systems).

Here, contrarily to  previous results obtained for this type of reaction-diffusion models, there are some regions in $(a,b)$ where $u$ is never
measured: we only require to know (i')
the density $u(0,x)$ in $(a,b)$ at $t=0$ and (ii') the density $u(t,x_0)$ and its spatial derivative $\ds{\frac{\partial u}{\partial x}}(t,x_0)$ for $t\in (0,\varepsilon)$ and some point $x_0$ in $(a,b)$ (see Remark  \ref{rem1} for a particular example of hypothesis (ii')). Thus a measurement of type (iii) is no more necessary. Furthermore, we show simultaneous uniqueness of two coefficients $\mu(x)$ and $\gamma$ provided that measurements of the second derivative $\ds{\frac{\partial^2 u}{\partial x^2}}(t,x_0)$ are available.

Our paper is organized as follows: in the next section, we give precise statements of our hypotheses and results; Section \ref{section_proof} is then dedicated to the proof of the results. Section \ref{section_num} is  devoted to the description of numerical examples illustrating how  the coefficient $\mu(x)$ can be approached using measures of the type (i') and (ii'). Those results are further discussed in Section \ref{sec:ccl}.

\section{Hypotheses and main results \label{sec:res}}

Let $(a,b)$ be an interval in $\R$. We consider the problem:
$$\left\{\baa{l}
\frac{\partial \ua}{\partial t}-D  \frac{\partial^2 \ua}{\partial x^2}=\ua\, (\mu(x)-\gamma \ua), \ t\ge 0, \ x\in(a,b),\\
\alpha_1 \ua (t,a) - \beta_1 \frac{\partial \ua}{\partial x} (t,a)=0, t>0, \\ \alpha_2 \ua (t,b) + \beta_2 \frac{\partial \ua}{\partial x} (t,b)=0, \ t>0,\\
\ua(0,x)=u_i(x), \ x\in(a,b).
\eaa
\right. \hspace{1cm} (\mathcal{P}_{\mu,\gamma})
$$
Our hypotheses on the coefficients are the following. Firstly, we assume that:
\be \label{hypmu}
\mu\in M:=\{\psi \in C^{0,\eta}([a,b]) \hbox{ such that }\psi \hbox{ is piecewise analytic on }(a,b)\},
\ee
for some $\eta \in (0,1]$. The space $C^{0,\eta}$ corresponds to H\"older continuous functions with exponent $\eta$ (see e.g. \cite{Fri64}). A function $\psi \in C^{0,\eta}([a,b])$ is called piecewise
analytic if it exists $n>0$ and an increasing sequence $(i_k)_{1\le k\le n}$ such that $i_1=a$, $i_n=b$, and $$\hbox{for all }x\in(a,b), \ \psi(x)=\sum_{j=1}^{n-1} \chi_{[i_j,i_{j+1})}(x)\varphi_j(x),$$for some analytic functions $\varphi_j$, defined on the intervals $[i_j,i_{j+1}]$, and where $\chi_{[i_j,i_{j+1})}$ are the characteristic functions of the intervals $[i_j,i_{j+1})$ for $j=1,\ldots, n-1$.

We also assume that $\gamma$ is a positive constant and that the boundary coefficients satisfy:
 \be \alpha_1,\alpha_2,\beta_1,\beta_2\ge 0 \hbox{ with }\alpha_1+\beta_1>0 \hbox{ and }\alpha_2+\beta_2>0. \label{HypBC}\ee

We furthermore make  the following hypotheses on the initial condition:

 \be u_i\ge 0, \ u_i\not \equiv 0 \hbox{ and }u_i\in C^{2,\eta}([a,b]),\label{HypDI}\ee
 for some $\eta$ in $(0,1)$, that is $u_i$ is a $C^2$ function such that $u_i''$ is H\"older continuous. In addition to that, we assume the following compatibility conditions:
\be\alpha_1 u_i (a) - \beta_1 u_i' (a)=0, \ \alpha_2 u_i (b) + \beta_2 u_i' (b)=0, \ \delta_{\beta_1}u_i''(a)=0, \ \delta_{\beta_2}u_i''(b)=0,\label{HypDI2}\ee
where $\delta_y$ is verifies: $\delta_0=1$ and $\delta_y=0$ if $y\neq 0$. We also need to assume that:
\be \label{hypui} \hbox{measure}(\{x \in (a,b), \ u_i(x)=0\})=0.
\ee

Under the assumptions (\ref{hypmu})-(\ref{HypDI2}), for each $\mu \in M$ and $\gamma>0$, the problem $(\mathcal{P}_{\mu,\gamma})$ has a unique solution $\ua\in C^{2,\eta}_{1,\eta/2}([0,+\infty)\times[a,b])$ (i.e. the derivatives up to order two in $x$ and order one in $t$ are H\"older continuous, see \cite{Fri64,Pao92} for a definition of H\"older continuity). Existence, uniqueness and regularity of the solution $u$ are classical. See e.g.  \cite[Ch. 1]{Pao92}.

\

Let us state our main results:
\begin{theorem}
\label{th:uniq1}
Let $\mu,\tmu \in M,$ and $\gamma>0$, and assume that the solutions $\ua$ and $\ub$ to $(\mathcal{P}_{\mu,\gamma})$ and $(\mathcal{P}_{\tmu,\gamma})$ satisfy, at some $x_0\in (a,b),$  and for some $\varepsilon>0$ and all  $t$ in $(0,\varepsilon)$: \ba  \ua(t,x_0) & = & \ub(t,x_0), \label{hypth1aa} \\ \frac{\partial \ua}{\partial x}(t,x_0) & = & \frac{\partial \ub}{\partial x} (t,x_0) . \label{hypth1ab} \ea Assume furthermore that
\be u_i(x_0)\neq 0 \hbox{ or } \frac{\partial^2 \ua}{\partial x^2}(t,x_0)=\frac{\partial^2 \ub}{\partial x^2} (t,x_0) \hbox{ for } t\in(0,\varepsilon). \label{hypth1b}\ee
 Then, we have $\mu\equiv\tmu$ on $[a,b]$ and consequently $\ua\equiv \ub$ in $[0,+\infty)\times[a,b]$. If $\beta_1>0$ (resp. $\beta_2>0$), this statement remains true when $x_0=a$ (resp. $x_0=b$).
\end{theorem}
\begin{remark}
This result remains valid if $\gamma=\gamma(x)$ is a given, positive function in $C^{0,\eta}([a,b]).$
\end{remark}
However, the conclusion of Theorem \ref{th:uniq1} is not true in general without the assumption  (\ref{hypth1ab}):
\begin{proposition}
\label{prop:uniq1}
Let $\mu \in M$ and $\gamma>0$. Assume that $\alpha_1=\alpha_2$  and $\beta_1=\beta_2$ and that $u_i$ is symmetric with respect to $x=(a+b)/2.$
Let $\tmu:=\mu(b-(x-a))$ for $x\in [a,b].$ Then, the solutions $\ua$ and $\ub$ to $(\mathcal{P}_{\mu,\gamma})$ and $(\mathcal{P}_{\tmu,\gamma})$ satisfy $u(t,\frac{a+b}{2}) = \ub(t,\frac{a+b}{2})$ for all $t\ge 0.$
\end{proposition}

Under an additional assumption on the initial condition $u_i$, we are able to obtain a uniqueness result for triples $(\ua,\mu,\gamma)$:
\begin{theorem}
\label{th:uniq2}
Let $\mu,\tmu \in M,$ and $\gamma,\tg>0$. Assume that, at some $x_0 \in (a,b)$, $u_i(x_0)=0$. Assume furthermore that the solutions $\ua$ and $\uc$ to $(\mathcal{P}_{\mu,\gamma})$ and $(\mathcal{P}_{\tmu,\tg})$ satisfy, for some $\varepsilon>0$ and for all  $t$ in $(0,\varepsilon)$: \ba \ua(t,x_0) & = &\uc(t,x_0), \label{hypth2a} \\ \frac{\partial \ua}{\partial x}(t,x_0)& = & \frac{\partial \uc}{\partial x} (t,x_0), \label{hypth2b} \\ \frac{\partial^2 \ua}{\partial x^2}(t,x_0)&= &\frac{\partial^2 \uc}{\partial x^2} (t,x_0). \label{hypth2c} \ea  Then, we have $\mu\equiv\tmu$ on $[a,b]$ and $\gamma=\tg$. Consequently $\ua\equiv \uc$ in $[0,+\infty)\times[a,b]$. If $\beta_1>0$ (resp. $\beta_2>0$), this statement remains true for $x_0=a$ (resp. $x_0=b$).
\end{theorem}

\begin{remarks}
\begin{itemize}
\item A particular example where hypotheses (\ref{hypth1aa}-\ref{hypth1b}) of Theorem \ref{th:uniq1} (resp. hypotheses (\ref{hypth2a}-\ref{hypth2c}) of Theorem \ref{th:uniq2}) are fulfilled is whenever, for some subset $\omega$ of $(a,b)$, $\ua(t,x)=\ub(t,x)$ for $t\in(0, \varepsilon)$ and all $x\in \omega$ (resp. $x_0\in \omega$ and $\ua(t,x)=\uc(t,x)$ in $(0,\varepsilon)\times \omega$). Note that, under this hypothesis, the previous results  \cite{CriRoq08} did not imply uniqueness; indeed, an additional assumption of type (iii) was required (cf. the introduction  section).

\item The uniqueness result of Theorem \ref{th:uniq2} cannot be adapted to the stationary equation associated to $(\mathcal{P}_{\mu,\gamma})$: $-p''=p\, (\mu(x)-\gamma p)$ (see e.g. \cite{BerHamRoq05a} for the existence and uniqueness of the stationary state $p>0$). Indeed, for any $\tau \in (0,1),$ setting $\tmu=\mu-\tau \gamma p$ and $\tg= (1-\tau)\gamma$, we obtain $-p''=p\, (\tmu(x)-\tg p)$, whereas $\tmu\not \equiv \mu$ and $\tg \not \equiv \gamma$. Thus, a measurement of $p$, even on the whole interval $[a,b]$, does not provide a unique couple $(\mu,\gamma)$.

\item The subset $M$ of $C^{0,\eta}([a,b])$ made of piecewise analytic functions is much larger than the set of analytic functions on $[a,b]$. It indeed contains some functions whose regularity is not higher than $C^{0,\eta}$, and some functions which are constant on some subsets of $[a,b]$. Our results hold true if  $M$ is replaced by any subset $M'$ of $C^{0,\eta}([a,b])$ such that for any couple of elements in $M'$, the subset of $[a,b]$ where these two elements intersect has a finite number of connected components.
\end{itemize}
\label{rem1}
\end{remarks}
\section{Proofs}\label{section_proof}

Let $\mu,\tmu \in M,$ and $\gamma,\tg>0$. Let $\ua$ be the solution to $(\mathcal{P}_{\mu,\gamma})$  and $\uc$ the solution to $(\mathcal{P}_{\tmu,\tg})$. We set $$U:=\ua-\uc \hbox{ and }m:=\mu-\tmu.$$ The function $U$ satisfies:

\begin{equation}
\frac{\partial U}{\partial t}- D  \frac{\partial^2 U}{\partial x^2}=\tmu U- \tg U(\ua+\uc)+\ua( m - \ua (\gamma-\tg)),
\label{eqU}
\end{equation}
for $t\ge 0$ and $x\in(a,b),$ and
\begin{equation}\left\{\baa{l}
\alpha_1 U (t,a) - \beta_1 \frac{\partial U}{\partial x} (t,a)=0, \ \alpha_2 U (t,b) + \beta_2 \frac{\partial U}{\partial x} (t,b)=0, \ t>0,\\
U(0,x)=0, \ x\in(a,b).
\eaa
\right.
\label{eqUb}
\end{equation}

\

\underline{Proof of Theorem \ref{th:uniq1}:} In that case $\gamma=\tg$. Equation (\ref{eqU}) then reduces to
\begin{equation}
\frac{\partial U}{\partial t}- D  \frac{\partial^2 U}{\partial x^2}=\tmu U- \gamma U(\ua+\ub)+\ua \, m.
\label{eqU1}
\end{equation}

\

\textit{Step 1: We prove that $m(x_0)=0.$}

\

It follows from hypothesis (\ref{hypth1aa}) that,  for all $t \in [0,\varepsilon)$, $U(t,x_0)=0$ and thereby, $$\frac{\partial U}{\partial t}(t,x_0)=0\hbox{ for all }t\in[0,\varepsilon).$$  If $u_i(x_0) \neq 0$, then, since $U(0,\cdot)\equiv 0$ we  deduce from (\ref{eqU1}) applied at $t=0$ and $x=x_0$ that $u_i(x_0) m(x_0)=0,$ and therefore $m(x_0)=0$.

If $u_i(x_0) =0$, from (\ref{hypth1b}), we have $\frac{\partial^2 U}{\partial x^2}(t,x_0)=0$ for all $t\in[0,\varepsilon)$. Applying equation (\ref{eqU1}) at $t=\varepsilon/2$ and $x=x_0,$ we get $$\ua\left(\frac{\varepsilon}{2},x_0\right) m(x_0)=0.$$ If $x_0 \in (a,b)$, the strong parabolic maximum principle (Corollary \ref{corPM}) applied to $\ua$ implies that $\ua(\varepsilon/2,x_0)>0$. As a consequence we again get $m(x_0)=0$. Lastly, if $x_0=a$ and $\beta_1>0$ the Hopf's Lemma applied to $\ua$ again implies that $\ua(\varepsilon/2,x_0)>0$. Indeed, assume on the contrary that $\ua(\varepsilon/2,x_0)=\ua(\varepsilon/2,a)=0$. The boundary condition $\alpha_1 \ua (t,a) - \beta_1 \frac{\partial \ua}{\partial x} (t,a)=0$ implies: $$\beta_1 \frac{\partial \ua}{\partial x} \left(\frac{\varepsilon}{2},a\right)=0,$$ which is impossible from Hopf's Lemma (Corollary \ref{corPM} and Theorem \ref{thPM} (b) and (c)). Thus $\ua(\varepsilon/2,x_0)>0$ and, again, $m(x_0)=0$. A similar argument holds for $x_0=b$, whenever $\beta_2>0.$

Under the assumptions of Theorem \ref{th:uniq1}, we therefore always obtain $m(x_0)=0.$

\

\textit{Step 2: We prove that $m\equiv 0$.}

\

Let us now set
$$b_1:=\sup\{x \in [x_0,b] \hbox{ s.t. }m\hbox{ has a constant sign on }[x_0,x]\}.$$By ``constant sign" we mean that either $m\ge0$ on $[x_0,x]$ or
$m\le 0$ on $[x_0,x]$. Then, four possibilities may arise:
\begin{itemize}
\item (i) $m=0$ on $[x_0,b_1]$ and $b_1<b$,\item (ii) $m\ge 0$ on $[x_0,b_1]$ and it exists $x_1 \in (x_0, b_1)$ such that $m(x_1)>0$,
\item (iii) $m\le 0$ on $[x_0,b_1]$ and it exists $x_1 \in (x_0, b_1)$ such that $m(x_1)<0$,
\item (iv) $b_1=b$, and $m=0$ on $[x_0,b]$.
\end{itemize}

Assume (i). Then, by definition of $b_1$, there exists a  decreasing sequence $y_k\to b_1,$ $y_k>b_1$, such that $|m(y_k)|>0$ for all $k  \ge 0$.
Assume that it exists $k_0$ such that $|m(x)|>0$ for all $x\in (b_1,y_{k_0})$. By continuity, $m$ does not change sign in $(b_1,y_{k_0})$, and therefore in $[x_0,y_{k_0}]$. This contradicts the definition of $b_1$. Thus, \be \label{zk}\hbox{ for all }k, \hbox{ it exists }z_k\in (b_1,y_k)\hbox{ such that }m(z_k)=0.\ee
Since $\mu$ and $\tmu$ belong to $M$, the function $m$ also belongs to $M$ and is therefore piecewise analytic on $(a,b)$. Thus, the set $\{x\in(a,b) \hbox{ s.t. }m(x)=0\}$ has a finite number of connected components. This contradicts (\ref{zk}) and rules out possibility (i).

Now assume (ii). By continuity of $m$, and from hypothesis (\ref{hypui}) on $u_i$, we can assume that $u_i(x_1)>0$. Since $m(x_1)>0$ and $U(0,\cdot)\equiv 0$, it follows from (\ref{eqU1})
 that $$\frac{\partial U}{\partial t}(0,x_1)= u_i(x_1) m(x_1)>0.$$ Thus, for $\varepsilon_1>0$ small enough, $U(t,x_1)>0$ for  $t\in(0,\varepsilon_1]$. As a consequence, $U$ satisfies:
\begin{equation}\left\{\baa{l}
\frac{\partial U}{\partial t}- D  \frac{\partial^2 U}{\partial x^2}-(\tmu-\gamma \ua- \gamma \uc) U \ge 0, \ t \in (0, \varepsilon_1], \ x\in (x_0,x_1),
\\
U (t,x_0)=0 \hbox{ and }U(t,x_1)>0 , \ t\in(0,\varepsilon_1], \\
U(0,x)=0, \ x\in (x_0,x_1).
\eaa
\right.
\label{eqU2}
\end{equation}
Moreover,
\begin{lemma}
We have $U(t,x)>0$ in $(0,\varepsilon_1)\times(x_0,x_1)$.
\label{lem1}
\end{lemma}
\textbf{Proof of Lemma \ref{lem1}:}
Set $W=U e^{-\lambda t},$ for some $\lambda>0$ large enough such that $c(t,x):=\tmu-\gamma \ua- \gamma \uc-\lambda \le 0$ in $(x_0,x_1).$ The function $W$ satisfies
$$\frac{\partial W}{\partial t}-D  \frac{\partial^2 W}{\partial x^2}-c(t,x) W \ge 0,  \ t \in (0, \varepsilon_1], \ x\in (x_0,x_1).$$
Assume that it exists a point $(t^*,x^*)$ in $(0,\varepsilon_1)\times (x_0,x_1)$ such that $U(t^*,x^*)< 0$. Then, since $W(t,x_0)=0$ and $W(t,x_1)>0$  for $t\in (0,\varepsilon_1)$, and since $W(0,x)=U(0,x)= 0,$ $W$ admits a minimum $m^*< 0$ in $(0,\varepsilon_1]\times (x_0,x_1)$. Theorem \ref{thPM} (a) applied to $W$ implies that $W\equiv m^*< 0$ on $[0,\varepsilon_1]\times [x_0,x_1]$, which is impossible. Thus $W\ge 0$
in $[0, \varepsilon_1]\times [x_0,x_1].$ Theorem \ref{thPM} (a) and (c) then implies that $W>0$ and consequently $U(t,x)>0$ in $(0,\varepsilon_1)\times(x_0,x_1)$. $\Box$

\

Since $U(t,x_0)=0$, the Hopf's lemma (Theorem \ref{thPM} (b) and (c)) also implies that $\ds{\frac{\partial U}{\partial x}(t,x_0)>0}$ for all $t\in (0,\varepsilon_1).$ This contradicts hypothesis (\ref{hypth1ab}). Possibility (ii) can therefore be ruled out.

Applying the same arguments to $-U$, possibility (iii) can also be rejected.
Finally, only (iv) remains.

Setting $$a_1:=\inf\{x\in [a,x_0] \hbox{ s.t. }m\hbox{ has a constant sign on }[x,x_0]\},$$ the same argument as above shows that $a_1=a$ and $m=0$ on $[a,x_0]$. Thus, finally, $m\equiv 0$ on $[a,b]$ and this concludes the proof of Theorem \ref{th:uniq1}. $\Box$

\

\underline{Proof of Theorem \ref{th:uniq2}:} From the assumptions (\ref{hypth2a}) and (\ref{hypth2c}) of Theorem \ref{th:uniq2},  equation  (\ref{eqU}) at $x=x_0$ reduces to
$$\ua(t,x_0)\, (m(x_0) - \ua(t,x_0) (\gamma-\tg))=0 \hbox{ for }t\in[0,\varepsilon).$$
If $x_0\in (a,b),$  the strong parabolic maximum principle (Corollary \ref{corPM}) implies that $\ua(t,x_0)>0$ for all $t>0$. This remains true if $x_0=a$ (if $\beta_1>0$) or $x_0=b$ (if $\beta_2>0$); cf. the proof of Theorem \ref{th:uniq1}. We therefore get:
\be
\label{eqU3}
m(x_0) = \ua(t,x_0) (\gamma-\tg) \hbox{ for }t\in(0,\varepsilon).
\ee
From the continuity of $t\mapsto \ua(t,x_0)$ up to $t=0$, we have $m(x_0)=u_i(x_0) (\gamma-\tg)$. Thus, $u_i(x_0)=0$ implies that  $m(x_0)=0$ which in turns implies from (\ref{eqU3}), and since $\ua(t,x_0)>0$ for $t>0$, that $\tg=\gamma$. The end of the proof is therefore similar to that of Theorem \ref{th:uniq1}. $\Box$

\begin{remark}
Extension of the arguments used in the previous proof to higher dimensions is not straightforward. Indeed, placing ourselves  in a bounded domain $\Omega$ of $\R^N,$ with $N\ge 2$, we may consider the largest region $\Omega_1$ in $\Omega$, containing $x_0$ and such that $m$ has a constant sign in $\Omega_1$. Consider in the above proof the possibility (ii)$_N$ (instead of (ii)): $m\ge 0$ on $\overline{\Omega_1}$ and it exists $x_1 \in \Omega$ such that $m(x_1)>0$.
    Then it exists a subset $\omega_1$ of $\Omega_1$, such that $x_0 \in \partial \omega_1$ and  $u(t,x)>0$ on a portion of $\partial \omega_1$. However, we cannot assert that $U(t,x)\ge 0$ on $\partial \omega_1$, and (ii)$_N$ can therefore not be ruled out as we did for (ii).
\end{remark}

\

\underline{Proof of Proposition \ref{prop:uniq1}:} Under the assumptions of Proposition \ref{prop:uniq1}, we observe that $\ub(t,b-(x-a))$ is a solution of $(\mathcal{P}_{\mu,\gamma}).$ In particular, by uniqueness, we have $$u(t,x)=\ub(t,b-(x-a)), \hbox{ for all }t\ge 0 \hbox{ and }x\in [a,b].$$It follows that $u(t,\frac{a+b}{2}) = \ub(t,\frac{a+b}{2})$ for all $t\ge 0.$ $\Box$

\section{Numerical computations}\label{section_num}

The purpose  of this section is to verify numerically that the measurements (\ref{hypth1aa}-\ref{hypth1ab}) of Theorem \ref{th:uniq1} allow to obtain a good approximation of the coefficient $\mu(x)$, when $\gamma$ is known.

Assuming that $\mu$ belongs to a finite-dimensional subspace $E\subset M$ and measuring  the distance between the measurements of
the solutions of $(\mathcal{P}_{\mu,\gamma})$ and $(\mathcal{P}_{\tmu,\gamma})$ through the function
$$G_\mu(\tmu)=\|\ua(\cdot ,x_0)-\ub(\cdot,x_0)\|_{L^2(0,\varepsilon)}+\|\frac{\partial \ua}{\partial x}(\cdot,x_0)-\frac{\partial \ub}{\partial x}(\cdot,x_0)\|_{L^2(0,\varepsilon)},$$  we look for  the coefficient $\mu(x)$ as a minimizer of the function $G_\mu$. Indeed, $G_\mu(\mu)=0$ and, from Theorem \ref{th:uniq1}, this is the unique global minimum of $G_\mu$ in $M$.

Solving $(\mathcal{P}_{\mu,\gamma})$ by a numerical method (see Appendix B) gives an approximate solution $u^h$. In our numerical
tests, we therefore replace $G_\mu$ by the discretized functional
$$\widehat{G}_\mu(\tmu):=\|\ua^h(\cdot ,x_0)-\ub^h(\cdot,x_0)\|_{L^2(0,\varepsilon)}+\|\frac{\partial \ua^h}{\partial x}(\cdot,x_0)-\frac{\partial \ub^h}{\partial x}(\cdot,x_0)\|_{L^2(0,\varepsilon)}.$$

\begin{remark}
Since $(\mathcal{P}_{\mu,\gamma})$ and $(\mathcal{P}_{\tmu,\gamma})$ are solved with the same (deterministic) numerical method, we have $\widehat{G}_\mu(\mu)=0$. Thus $\mu$ is a global minimizer of $\widehat{G}_\mu$. However, this minimizer  might not be unique.
\end{remark}

\subsection{State space $E$ \label{subsec_E}}

We fix $(a,b)=(0,1)$ and we assume that the function $\mu$ belongs to a subspace  $E \subset M$ defined by:
$$E:=\left\{\tmu \in  C^{0,\eta}([0,1]) \, | \, \exists \  (h_i)_{0\le i \le n}\in \R^{n+1}, \ \tmu(x)=\sum_{i=0}^{n}h_i \cdot j\left((n-2)\left(x-c_i\right)\right) \hbox{ on }[0,1]\right\},$$ with $c_i=\frac{i-1}{n-2}$ and
$j(x)=\left\{\baa{l} \exp\left(\frac{4 x^2}{x^2-4}\right), \hbox{ if }x\in(-2,2), \\ 0 \hbox{ otherwise.}
\eaa \right.$

\subsection{Minimization of $\widehat{G}_\mu$ in $E$ \label{subsec_minG}}

For the numerical computations,  we fixed  $D=0.1$, $\gamma=1$, $\alpha_1=\alpha_2=0$ and $\beta_1=\beta_2=1$ (Neumann boundary conditions). Besides, we assumed that $u_i\equiv 0.2$, $\varepsilon=0.3$ and $x_0=2/3.$ The integer $n$ was set to $10$ in the definition of $E$.

Numerical computations were carried out for $100$ functions $\mu_k$ in $E:$ $$\mu_k=\sum_{i=0}^{n}h_i^k \cdot j\left[(n-2)\left(x-c_i\right)\right], \ k=1\ldots 100,$$whose components $h_i^k$ were randomly drawn from  a uniform distribution in $(-5,5).$

Minimizations of the functions $\widehat{G}_{\mu_k}$ were performed using MATLAB's$^{\circledR}$ \textit{fminunc} solver \footnote{ MATLAB's$^{\circledR}$ \textit{fminunc} medium-scale optimization algorithm uses a Quasi-Newton method with a mixed quadratic and cubic line search procedure. Our stopping criterion was based on the maximum number of evaluations of the function $\widehat{G}_\mu$, which was set at $2\cdot 10^3.$}. This led to $100$ functions $\mu^*_k$ in $E$, each one corresponding to a
computed approximation for a minimizer of
the function $\widehat{G}_{\mu_k}$. In our numerical tests, we obtained values of $\widehat{G}_{\mu_k}(\mu^*_k)$ in $(5\cdot 10^{-7}, 10^{-5})$, with an average of $5\cdot 10^{-6}$ and a standard deviation of $2\cdot 10^{-6}.$

The values $\|\mu_k-\mu^*_k\|_{L^2([0,1])}/\|\mu_k\|_{L^2([0,1])}$, for $k=1\ldots 100$, are comprised between $5 \cdot 10^{-3}$ and $0.16$, with an average value of $0.04$ and a standard deviation of $0.03$.

%This average distance is low compared to the error made by approaching $\mu_k$ by the constant value $\mu_k(x_0).$ Indeed, the average value of $\|\mu_k-\mu_k(x_0)\|_{L^2([0,1])}$ is $3.45$.

\begin{figure}
\centering
\subfigure[$\mu$ (plain line) and $\mu^*$ (dotted line)]{
\includegraphics*[width=7cm]{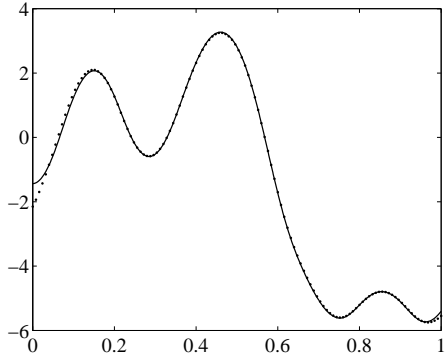}}
\subfigure[$\mu$ (plain line) and $\underline{\mu}^*$ (dotted line)]{ \includegraphics*[width=7cm]{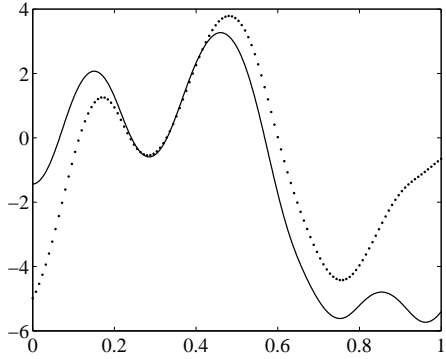}}
\caption{(a) An example of function $\mu$ in $E$, together with a function $\mu^*$ which was obtained by minimizing $\widehat{G}_{\mu}.$ In this case  $\|\mu-\mu^*\|_{L^2([0,1])}/\|\mu\|_{L^2([0,1])}=0.03$ and $\widehat{G}_{\mu}(\mu^*)=2 \cdot 10^{-6}$. (b) The same function $\mu$ together with the function $\underline{\mu}^*$ obtained by minimizing  $\widehat{H}_{\mu}$. Here, $\|\mu-\underline{\mu}^*\|_{L^2([0,1])}/\|\mu\|_{L^2([0,1])}=0.47$ and $\widehat{H}_{\mu}(\underline{\mu}^*)=2 \cdot 10^{-6}.$}
\label{fig:exple}
\end{figure}

Fig. \ref{fig:exple} (a) depicts an example of function $\mu$ in $E$, together with a function $\mu^*$ which was obtained by minimizing $\widehat{G}_{\mu}.$

\

\subsection{Test of another criterion $H_\mu$ \label{subsec_minH}}

In this section, we illustrate that measurement (\ref{hypth1aa}) alone  cannot be used for reconstructing $\mu$. Replacing $G_\mu$ by:
$$H_\mu(\tmu)=\|\ua(\cdot,x_0)-\ub(\cdot,x_0)\|_{L^2(0,\varepsilon)},$$and setting
$\widehat{H}_\mu(\tmu):=\|\ua^h(\cdot ,x_0)-\ub^h(\cdot,x_0)\|_{L^2(0,\varepsilon)},$
we performed the same analysis as above, with the same samples $\mu_k\in E$ and the same parameters.

The corresponding values of $\widehat{H}_{\mu}(\underline{\mu}^*_k)$ are comparable to those obtained in Section \ref{subsec_minG}.  Namely, these values are included in $(2\cdot 10^{-8}, 10^{-5}),$ with average $2\cdot 10^{-6},$ and standard deviation $3\cdot 10^{-6}.$
However, the corresponding values of the distance $\|\mu_k-\underline{\mu}^*_k\|_{L^2([0,1])}/\|\mu_k\|_{L^2([0,1])}$ are far larger than those obtained in Section \ref{subsec_minG}: these values are comprised between $0.08$ and $1.64$, with an average of $0.56$ and a standard deviation of $0.34$.

Using the same sample $\mu \in E$ as in Fig. \ref{fig:exple} (a), we present in Fig. \ref{fig:exple} (b) the approximation $\underline{\mu}^*$ obtained by minimizing $\widehat{H}_{\mu}$. In this case, the distance $\|\mu-\underline{\mu}^*\|_{L^2([0,1])}/\|\mu\|_{L^2([0,1])}$ is $18$ times larger than $\|\mu-\mu^*\|_{L^2([0,1])}/\|\mu\|_{L^2([0,1])}.$

\section{Discussion \label{sec:ccl}}

Studying the reaction-diffusion problem $(\mathcal{P}_{\mu,\gamma})$ with a nonlinear term of the type $u \, (\mu(x)-\gamma u),$ we have
proved in Section \ref{sec:res} that  knowing  $u$ and its first spatial derivative
at a single point $x_0$ and for small times $t\in (0,\varepsilon)$ is sufficient to completely determine the couple $(u(t,x),\mu(x))$. Additionally, if the second spatial derivative is also measured at $x_0$ for $t\in (0,\varepsilon)$, the triplet $(u(t,x),\mu(x),\gamma)$ is also completely determined.

These uniqueness results  are mainly the consequences of Hopf's Lemma and of an hypothesis on the set $M$ of coefficients which $\mu(x)$ belongs to. This hypothesis implies that two coefficients in $M$ can be equal only over a set having a finite number of connected components.

The theoretical results of Section \ref{sec:res} suggest that the coefficients $\mu(x)$ and $\gamma$ can be numerically determined using only measurements of the solution $u$ of $(\mathcal{P}_{\mu,\gamma})$ and of its spatial derivatives at one point $x_0,$ and for $t\in (0,\varepsilon)$.
Indeed, the numerical computations of Section \ref{section_num} show that, when $\gamma$ is known,
the coefficient $\mu(x)$ can be estimated  by minimizing a function $G_\mu$. The function $\ub$ being the solution of $(\mathcal{P}_{\tmu,\gamma}),$ we defined  $G_\mu(\tmu)$ as the distance between  $(u,\partial u/ \partial x)(\cdot,x_0)$ and $(\ub,\partial \ub/\partial x)(\cdot,x_0),$ in the $L^2(0,\varepsilon)$ sense.

The numerical computations presented in Section \ref{subsec_minG} were carried out on $100$ samples of functions $\mu_k$ chosen in a finite-dimensional subspace of $M$. In each case, a good approximation $\mu^*_k$ of $\mu_k$ was obtained. The average relative $L^2$-error  between $\mu^k$ and $\mu_k^*$ is $30$ times smaller than the average relative $L^2$-error between $\mu^k$ and the constant function $\mu_k(x_0)$.
Thus, a measurement of $u$ and of its first spatial derivative at a point $x_0$ (and for $t\in (0,\varepsilon)$) indirectly gives more information on the global shape of $\mu$ than a direct measure of $\mu$ at $x_0.$
These good results, in spite of the computational error, indicate $L^2$-stability of the coefficient $\mu$ with respect to single-point measurements of the solution $u$ of $(\mathcal{P}_{\mu,\gamma})$ and of its spatial derivative.

Proposition \ref{prop:uniq1} shows that the uniqueness result
of Theorem \ref{th:uniq1} is not true
without the assumption (\ref{hypth1ab}) on the spatial derivatives. This suggests that measurement (\ref{hypth1aa})  alone  cannot be used for reconstructing $\mu.$ In Section
\ref{subsec_minH}, working with the same samples $\mu_k$ as those discussed
above, we obtained approximations $\underline{\mu}^*_k$ of
$\mu_k$ by minimizing a new function $H_\mu$, which measures  the distance between $u(\cdot,x_0)$ and $\ub(\cdot,x_0).$
The  average relative $L^2$-error between $\mu_k$ and $\underline{\mu_k}^*$  was $14$ times larger than the average relative $L^2$-error separating   $\mu_k$ and $\mu^*_k.$ This confirms the usefulness of the spatial derivative measurements for the reconstruction of $\mu.$

\section*{Acknowledgements}

The authors would like to thank two anonymous referees for their valuable comments on an earlier version of this paper.
The first author is supported by the French ``Agence Nationale de la Recherche" within the projects ColonSGS, PREFERED, and URTICLIM.

\section*{Appendix A: maximum principle}

The following version of the parabolic maximum principle can be found in \cite[Ch. 2]{Fri64} and \cite[Ch. 3]{ProWei67}.

\begin{theorem}
Let $u \in C^{2}_{1}((0,T]\times (x_1,x_2))\cap C([0,T]\times [x_1,x_2]),$ for some $T>0$ and $x_1,x_2\in\R.$ Let $c(t,x)\le 0 \in C^{0,\eta}_{0,\eta/2}([0,T]\times [x_1,x_2])$, for some $\eta \in (0,1]$.

Suppose that $\ds{\frac{\partial u}{\partial t}-D  \frac{\partial^2 u}{\partial x^2}-c(x) u \ge 0}$ for  $t\in(0,T]$ and $x\in (x_1,x_2)$.

\

(a) If $u$ attains a minimum $m^* \le 0$ at a point $(t^*,x^*)\in (0,T]\times (x_1,x_2),$ then $u(t,x)\equiv m^*$ on $[0,t^*]\times [x_1,x_2]$.

(b) (Hopf's Lemma) If $u$ attains a minimum $m^* \le 0$ at a point $(t^*,x_1)$ (resp. $(t^*,x_2)$), with $t^*>0$, then either $\frac{\partial u}{\partial x}(t^*,x_1)>0$ (resp. $\frac{\partial u}{\partial x}(t^*,x_2)<0$) or $u(t,x)\equiv m^*$ on $[0,t^*]\times [x_1,x_2]$.

(c) If $u\ge 0,$ the results (a) and (b) remain true without the assumption $c(t,x)\le 0.$
\label{thPM}
\end{theorem}

\noindent An immediate corollary of this theorem is:
\begin{corollary}
The solution $u(t,x)$ of $(\mathcal{P}_{\mu,\gamma})$ is strictly positive in $(0,+\infty)\times(a,b)$.
\label{corPM}
\end{corollary}

\textbf{Proof of Corollary \ref{corPM}:} Assume that it exists $(t^*,x^*)\in (0,+\infty)\times (a,b)$ such that $u(t^*,x^*)<0$.

Set $w(t,x)=u  \, e^{-\lambda t},$ for $\lambda>0$ large enough such that $$c(t,x):=\mu(x) - \gamma u - \lambda\le 0 \hbox{ in } [0,t^*]\times [a,b].$$The function $w$ satisfies:
$$\frac{\partial w}{\partial t}-D  \frac{\partial^2 w}{\partial x^2}-c(t,x) w =0.$$Since $w(0,x)=u_i(x)\ge 0$ in $(a,b)$ and $w(t^*,x^*)<0,$ the function $w$ admits a minimum $m^*<0$ in $(0,t^*]\times [a,b].$ From Theorem \ref{thPM} (a), and since $u_i\not \equiv 0,$ this minimum is attained at a boundary point: it exits $t'\in (0,t^*]$ such that $w(t',a)=m^*<0$ or $w(t',b)=m^*<0$. Without loss of generality, we can assume in the sequel that $w(t',a)=m^*<0.$ From Theorem \ref{thPM} (b), we obtain $\frac{\partial w}{\partial x}(t',a)>0.$ Using the boundary conditions in problem $(\mathcal{P}_{\mu,\gamma})$, we finally get: $$\alpha_1 \, m^* = \beta_1 \, \frac{\partial w}{\partial x}(t',a)>0.$$Using assumption (\ref{HypBC}), we get a contradiction. Thus $u(t,x)\ge 0$ in $(0,+\infty)\times (a,b).$ The conclusion then follows from Theorem \ref{thPM} (c). $\Box$

\

\section*{Appendix B: numerical solutions of $(\mathcal{P}_{\mu,\gamma})$ and $(\mathcal{P}_{\tmu,\gamma})$}

The equations $(\mathcal{P}_{\mu,\gamma})$ and $(\mathcal{P}_{\tmu,\gamma})$ were solved using Comsol Multiphysics$^{\circledR}$ time-dependent solver, using second order finite element method (FEM) with 960 elements. This solver uses a method of lines approach incorporating variable order variable stepsize backward differentiation formulas. Nonlinearities are treated using a Newton's method. The interested reader can get more information in Comsol Multiphysics$^{\circledR}$ user's guide.

%\bibliographystyle{unsrt}
%\bibliography{biblio_lionel}

\bigskip

\section*{References}

\begin{thebibliography}{10}

\bibitem{Hol93}
E~E Holmes.
\newblock {Are diffusion models too simple? A comparison with telegraph models
  of invasion}.
\newblock {\em {American Naturalist}}, {142}:{779--795}, {1993}.

\bibitem{ShiKaw97}
N~Shigesada and K~Kawasaki.
\newblock {\em {Biological invasions: theory and practice}}.
\newblock {Oxford Series in Ecology and Evolution, Oxford: Oxford University
  Press}, {1997}.

\bibitem{Tur98}
P~Turchin.
\newblock {\em Quantitative analysis of movement: measuring and modeling
  population redistribution in animals and plants}.
\newblock {Sinauer Associates, Sunderland, MA}, {1998}.

\bibitem{Mur02}
J~D Murray.
\newblock {\em {Mathematical Biology}}.
\newblock {Third Edition. Interdisciplinary Applied Mathematics 17,
  Springer-Verlag, New York}, {2002}.

\bibitem{OkuLev02}
A~Okubo and S~A Levin.
\newblock {\em Diffusion and ecological problems -- modern perspectives}.
\newblock {Second edition, Springer-Verlag, New York}, {2002}.

\bibitem{Ske51}
J~G Skellam.
\newblock {Random dispersal in theoretical populations}.
\newblock {\em {Biometrika}}, {38}:{196--218}, {1951}.

\bibitem{PeaBla06}
K~Pearson and J~Blakeman.
\newblock {\em {Mathematical contributions of the theory of evolution. A
  mathematical theory of random migration}}.
\newblock Drapersi Company Research Mem. Biometrics Series~III, Dept. Appl.
  Meth. Univ. College, London, {1906}.

\bibitem{Fis37}
R~A Fisher.
\newblock {The wave of advance of advantageous genes}.
\newblock {\em {Annals of Eugenics}}, {7}:{335--369}, {1937}.

\bibitem{kpp}
A~N Kolmogorov, I~G Petrovsky, and N~S Piskunov.
\newblock {\'Etude de l'\'equation de la diffusion avec croissance de la
  quantit\'e de mati\`ere et son application \`a un probl\`eme biologique}.
\newblock {\em {Bulletin de l'Universit\'e d'\'Etat de Moscou, S\'erie
  Internationale A}}, {1}:{1--26}, {1937}.

\bibitem{ShiKaw86}
N~Shigesada, K~Kawasaki, and E~Teramoto.
\newblock {Traveling periodic-waves in heterogeneous environments}.
\newblock {\em {Theoretical Population Biology}}, {30}({1}):{143--160}, {1986}.

\bibitem{BerHamRoq05a}
H~Berestycki, F~Hamel, and L~Roques.
\newblock {Analysis of the periodically fragmented environment model: I -
  Species persistence}.
\newblock {\em Journal of Mathematical Biology}, 51(1):75--113, 2005.

\bibitem{CanCos03}
S~Cantrell, R and C~Cosner.
\newblock {\em {Spatial ecology via reaction-diffusion equations}}.
\newblock {John Wiley \& Sons Ltd, Chichester, UK }, {2003}.

\bibitem{ElsHamRoq09}
M~El~Smaily, F~Hamel, and L~Roques.
\newblock Homogenization and influence of fragmentation in a biological
  invasion model.
\newblock {\em {Discrete and Continuous Dynamical Systems, Series A}},
  25:321--342, 2009.

\bibitem{HamFay10}
F~Hamel, J~Fayard, and L~Roques.
\newblock Spreading speeds in slowly oscillating environments.
\newblock {\em {Bulletin of Mathematical Biology}}, {DOI
  10.1007/s11538-009-9486-7}, 2010.

\bibitem{RoqHam07}
L~Roques and F~Hamel.
\newblock Mathematical analysis of the optimal habitat configurations for
  species persistence.
\newblock {\em Mathematical Biosciences}, 210(1):34--59, 2007.

\bibitem{RoqSto07}
L~Roques and R~S Stoica.
\newblock Species persistence decreases with habitat fragmentation: an analysis
  in periodic stochastic environments.
\newblock {\em Journal of Mathematical Biology}, 55(2):189--205, 2007.

\bibitem{Xin00}
J~Xin.
\newblock {Front propagation in heterogeneous media}.
\newblock {\em {SIAM Review}}, {42}:{161--230}, {2000}.

\bibitem{BukKli81}
A~L Bukhgeim and M~V Klibanov.
\newblock {Global uniqueness of a class of multidimensional inverse problems}.
\newblock {\em {Soviet Mathematics - Doklady}}, {24}:{244--247}, {1981}.

\bibitem{Isa93}
V~Isakov.
\newblock {An uniqueness in inverse problems for semilinear parabolic
  equations}.
\newblock {\em {Archive for Rational Mechanics and Analysis}}, {1993}.

\bibitem{KliTim04}
M~V Klibanov and A~Timonov.
\newblock {\em {Carleman estimates for coefficient inverse problems and
  numerical applications}}.
\newblock {Inverse And Ill-Posed Series, VSP, Utrecht}, {2004}.

\bibitem{CriRoq08}
M~Cristofol and L~Roques.
\newblock Biological invasions: Deriving the regions at risk from partial
  measurements.
\newblock {\em Mathematical Biosciences}, 215(2):158--166, 2008.

\bibitem{ImmYam98}
O~Y Immanuvilov and M~Yamamoto.
\newblock {Lipschitz stability in inverse parabolic problems by the Carleman
  estimate}.
\newblock {\em {Inverse Problems}}, {14}:{1229--1245}, {1998}.

\bibitem{YamZou01}
M~Yamamoto and J~Zou.
\newblock {Simultaneous reconstruction of the initial temperature and heat
  radiative coefficient}.
\newblock {\em {Inverse Problems}}, {17}:{1181--1202}, {2001}.

\bibitem{BelYam06}
M~Belassoued and M~Yamamoto.
\newblock {Inverse source problem for a transmission problem for a parabolic
  equation}.
\newblock {\em {Journal of Inverse and Ill-Posed Problems}},
  {14}({1}):{47--56}, {2006}.

\bibitem{CriGai06}
M~Cristofol, P~Gaitan, and H~Ramoul.
\newblock {Inverse problems for a two by two reaction-diffusion system using a
  carleman estimate with one observation}.
\newblock {\em {Inverse Problems}}, {22}:{1561--1573}, {2006}.

\bibitem{BenCri08}
A~Benabdallah, M~Cristofol, P~Gaitan, and M~Yamamoto.
\newblock {Inverse problem for a parabolic system with two components by
  measurements of one component}.
\newblock {\em {Applicable Analysis}}, {88}(5):{683--710}, {2009}.

\bibitem{Fri64}
A~Friedman.
\newblock {\em { Partial differential equations of parabolic type}}.
\newblock {Prentice-Hall, Englewood Cliffs, NJ}, {1964}.

\bibitem{Pao92}
C~V Pao.
\newblock {\em {Nonlinear Parabolic and Elliptic Equations}}.
\newblock {Plenum Press, New York}, {1992}.

\bibitem{ProWei67}
M~H Protter and H~F Weinberger.
\newblock {\em {Maximum Principles in Differential Equations}}.
\newblock {Prentice-Hall, Englewood Cliffs, NJ}, {1967}.

\end{thebibliography}

\end{document}